\documentclass{math_crypt_v06}

\usepackage{amsfonts}

\title{Statistics of
Random Permutations and the
Cryptanalysis Of Periodic Block Ciphers}
\headlinetitle{Statistics of Random Permutations and Block Ciphers}

\authorone{Nicolas T. Courtois}
\addressone{University College of London, Gower Street, London, WC1E6BT}
\countryone{UK}
\emailone{n.courtois@ucl.ac.uk}
 
\authortwo{Gregory V. Bard}
\addresstwo{Fordham University, Department of Mathematics, The Bronx, NY, 10458}
\countrytwo{USA}
\emailtwo{bard@fordham.edu}

\authorthree{Shaun V. Ault}
\addressthree{Fordham University, Department of Mathematics, The Bronx, NY, 10458}
\countrythree{USA}
\emailthree{ault@fordham.edu}

\abstract{A block cipher is intended to be computationally indistinguishable from a 
random permutation of appropriate domain and range. But what are the properties
of a random permutation? By the aid of exponential and ordinary generating functions,
we derive a series of collolaries of interest to the cryptographic community. These
follow from the Strong Cycle Structure Theorem of permutations, and are useful in
rendering rigorous two attacks on Keeloq, a block cipher in wide-spread use. These
attacks formerly had heuristic approximations of their probability of success.

Moreover, we delineate an attack against the (roughly) millionth-fold iteration
of a random permutation. In particular, we create a distinguishing attack, whereby
the iteration of a cipher a number of times equal to a particularly chosen highly-composite 
number is breakable, but merely one fewer round is considerably \emph{more} secure.
We then extend this to a key-recovery attack in a ``Triple-DES'' style construction,
but using AES-256 and iterating the middle cipher (roughly) a million-fold. 

It is hoped that these results will showcase the utility of exponential and ordinary
generating functions and will encourage their use in cryptanalytic research.}
\acknowledgments{We thank Sebastiaan Indesteege, a graduate student
from {Katho-lieke} Universiteit Leuven in Belgium, for helpful comments;  
Sean O'Neil, an independent scientist from Ireland, also for
helpful comments; Dr. Kenneth Patterson of the 
Royal Holloway University of London for questions at the 2008 
Workshop on Mathematical Cryptography in Santander, Spain,
that encouraged us to rigorously describe the attack in 
Section~\ref{highlyiterated}; and Prof. Philippe Flajolet and Prof. 
Robert Sedgewick, for their excellent text on analytic combinatorics
\cite{Flajolet}.}

\firstpage{1}                                %
\doi{xxx}                                    %
\communicated{xxx}                           %
\received{xxx}                               %
\revised{xxx}                                %

\begin{document}




\newcommand{\integer}{\mathbb{Z}}
\newcommand{\fancy}[1]{{\cal #1}}
\newcommand{\set}[1]{\left \{ #1 \right \}}
\newcommand{\ignore}[1]{}
\newcommand{\floor}[1]{\left \lfloor #1 \right \rfloor}

\newcommand{\defequal}{\stackrel{def}{=}}


\ignore{
\def\N{\mbox{\rm{I}\hspace{-.15em}\rm{N}}}
\def\Z{\mbox{\rm{Z}\hspace{-.28em}\rm{Z}}}
\def\RR{\mbox{\rm{I}\hspace{-.15em}\rm{R}}}
\def\C{\mbox{\rm{l}\hspace{-.50em}\rm{C}}}
\def\Q{\mbox{\rm{l}\hspace{-.50em}\rm{Q}}}
\def\F{\mbox{\rm{I}\hspace{-.20em}\rm{F}}}

\spnewtheorem{defi}{Definition}[section]{\bfseries}{\rmfamily}
\spnewtheorem{prop}{Proposition}[section]{\bfseries}{\rmfamily}
\spnewtheorem{fact}{Fact}[section]{\bfseries}{\rmfamily}
\spnewtheorem{coro}{Corollary}[prop]{\bfseries}{\rmfamily}
}
\section{Introduction}
The technique of using a function of a variable to count objects of
various sizes, using the properties of multiplication and addition of
series as an aid, is accredited to Pierre-Simon Laplace \cite{Flajolet}. 
Here, we will use this family of techniques, now called ``analytic combinatorics''
to count permutations of particular types.
An ordinary
generating series associated with a set of objects assigns as the 
coefficient of the $z^i$th term, the number of objects of size $i$.
An exponential generating series is merely this, with each term
divided by $i!$. In particular, this can be used to describe permutations
drawn at random from $S_n$, which is the topic of this paper.

The cipher Keeloq, can be written as the eighth iterate of a 
permutation followed by one more permutation \cite[Ch. 2]{BardDis}. 
This eighth power naturally affects the cycle structure; for
example, we will prove that the fixed points of the eighth power
are those of order $\set{1,2,4,8}$ under the original. There are
many other properties of these repeated permutations that
follow from the factorization of the number of iterations, and
we will show cryptanalytic consequences.

In the remainder of this section we will introduce analytic
combinatorics through exponential and ordinary generating
functions. In Section~2 we prove a theorem on the cycle
structures of random permutations, and in Section~3 we
present a number of corollaries. We imagine that most if
not all of these are already known in some form, 
but here we are compiling
them all in one place, with a view to determining when a random
permutation has a given property, rather than merely counting
objects which is the usual use of the techniques of this subject.
The proofs are our own. In Section~4, we apply these techniques
to Keeloq, and describe two quite feasible attacks, but also their
\emph{exact} success probabilities. These attacks have been 
previously described as requiring the entire code-book of the
cipher (all plaintext-ciphertext pairs under the current key) but
here we let $\eta$ represent the fraction of the code-book available,
and show how $\eta$ affects the success probability. In Section~5,
we present an unusual example, where a very highly iterated cipher
appears to be secure, but adding one iteration opens up a feasible
and effective distinguishing attack. We conclude in Section~6.

\ignore{
\documentclass[11pt]{article}

\usepackage{amsfonts}
\usepackage{graphicx}
\usepackage[text={6.5in,9.25in},foot=0.6in]{geometry}
\usepackage{pstricks}
\usepackage{pst-plot}
\usepackage{pst-node}


\setlength{\parskip}{1.0ex plus0.2ex minus0.2ex}
\setlength{\parindent}{0.0in}
\renewcommand{\baselinestretch}{1.2}


\newtheorem{theorem}{Theorem \rm}
\newtheorem{lemma}{Lemma}
\newenvironment{proof}{\bf{Proof:}\rm}{$\square$\newline\newline}


\newcommand{\R}{\mathbb{R}}
\newcommand{\ds}{\displaystyle}


\begin{document}

{\Large 1.1 Background}
}

\newcommand{\Z}{\mathbb{Z}}

\subsection{Background}
A combinatorial class $\fancy{C}$ is a set of objects C together with a function 
$\ell_C : C \to \Z^{\geq 0}$, which asssigns to each element a non-negative integer ``size''. For 
example, if $P$ is the set of permutation groups $S_n$ for all positive integers n, then we may
use the size function $\ell_P(\pi) = n$, for any $\pi \in S_n$, to make $\fancy{P}$
into a combinatorial class.

Let $C_i$ be the cardinality of the set of elements of $C$ with size $i$. Thus in our example,
$P_i = i!$ for $i \geq 0$. It will be useful to represent $C_i$ by either an exponential or
an ordinary generating 
function (OGF or EGF).  First, a brief discussion of generating functions is in order.

Given a set of constants indexed by $\Z^{\geq 0}$, say $c_0, c_1, c_2, \ldots$, the {\it
ordinary generating function} (or OGF) is defined as the formal power series:
\[
  c(z) \defequal \sum_{i=0}^{\infty} c_iz^i = c_0 + c_1 z + c_2 z^2 + c_3 z^3 + \cdots.
\]
The EGF is defined as the formal power series:
\[
  c_e(z) \defequal \sum_{i=0}^{\infty} \frac{c_i}{i!} z^i = c_0 + \frac{c_1}{1!} z + 
  \frac{c_2}{2!} z^2 + \frac{c_3}{3!} z^3 + \cdots.
\]
For our example combinatorial class, $\fancy{P}$, its OGF is $\fancy{P}(z) = z + 2z^2 + 6z^3 + 24z^4 + 120z^5
+ \cdots$, and its EGF is $\fancy{P}_e(z) = z + z^2 + z^3 + z^4 + z^5 + \cdots$.  The
series $1 + z + z^2 + z^3 + z^4 + z^5 + \cdots$ represents the OGF of the 
non-negative integers, $\Z^{\ge 0}$ with
``size'' function being the identity: $\ell(n) = n$.

In combinatorial arguments,
OGFs and EGFs abound \cite{Flajolet} \cite{MarkoRiedel} and are
especially useful in counting partitions of sets.  For example, let $A_1$, $A_2$, \ldots, $A_k$
be sets of whole numbers.  The number of all distinct ways
that $n$ identical objects can be placed into 
$k$ containers, where container $j$ must have some number of objects that occurs in the set $A_j$
will be the coefficient of $z^n$ in the OGF:\label{statedhere}
\[
  \left( \sum_{i \in A_1} z^i \right)\left( \sum_{i \in A_2} z^i \right) \cdots
  \left( \sum_{i \in A_k} z^i \right),
\]
a function that we will use in the proof of Lemma~\ref{lem.B_is_k}.
Notice that the $j^{th}$ factor is the OGF that represents the set $A_j$.
There is a similar interpretation for EGFs and products of EGFs, in
terms of probability rather than strict counting. See Section~\ref{limcyclecount}
or Theorem~\ref{tautheorem} as an example.

A less trivial example of a combinatorial class is the class $\fancy{O}$ 
of $n$-cycles of $S_n$, for all $n > 0$, with size function
$\ell(\pi) = n$ if $\pi \in S_n$.  In other words, size $n$ members of $\fancy{O}$ comprise the subset
of permutations of $S_n$ where the permutation has exactly one orbit. For any $n > 0$ there are
$n!/n$ or $(n-1)!$ of these.  Thus the OGF is $z + z^2 + 2z^3 + 6z^4 + 24z^5 + 120z^6 + \cdots$ ,
and the EGF is $z + z^2/2 + z^3/3 + z^4/4 + z^5/5 + z^6/6 + \cdots$.
Thus the probability that a random permutation from $S_n$ has only one cycle is given
by the coefficients of the $z^n$ terms in the EGF.  Namely, $(n-1)!/n! = 1/n$.

Often, the formal power series defining OGFs or EGFs converge to functions (in some neighborhood
of $0$).  For example, the OGF for $\Z^{\geq 0}$ converges to $1/(1-z)$, and its EGF converges
to $e^z$.  The EGF for the combinatorial class $\fancy{O}$ above
\ignore{whose coefficients compute the probability that a randomly chosen 
permutation of $n$ letters is an $n$-cycle} also converges:
\[
  z + \frac{z^2}{2} + \frac{z^3}{3} + \frac{z^4}{4} + \frac{z^5}{5} + \frac{z^6}{6} + \cdots 
  = \log\left( \frac{1}{1-z} \right),
\]
as can be verified by term-by-term integration of the power series for $\frac{1}{1-z}$.  The 
existence of such functions will facilitate multiplications and compositions.

\subsection{Notation}

The somewhat unusual notation of $exp(C)$ where $C$ is a series, means precisely substituting
the entire series $C$ for $z$ into the Taylor expansion for $e^z = \sum_{i \geq 0} z^i/i!$, 
similar to matrix exponentiation. 

It is well-known that any permutation may be written uniquely as a product of 
disjoint cycles, up to reordering of the cycles and cyclic reordering within each cycle; 
indeed, for any given permutation $\pi$ consisting of $k$ disjoint cycles, 
having cycle lengths $c_1$, $c_2, c_3, \ldots,
c_k$, there are exactly $k! c_1c_2c_3\cdots c_k$ ways to reorder to obtain an equivalent
expression for $\pi$.  Any counts we make of symmetric group elements must take this fact into
account.  Note, we use the convention that if $\pi$ has a fixed-point, $a$, then the $1$-cycle
$(a)$ is part of the expression for $\pi$ as disjoint cycles.  In particular, the identity of
$S_n$ is written $(1)(2)(3)\cdots (n)$.  We use the term {\it cycle-count} for the number of 
disjoint cycles (including all $1$-cycles) in the expression of a permutation.
It shall be convenient to include in our analysis the unique permutation of no
letters, which has by convention
cycle-count $0$.  We may view this element as the sole member of $S_0$.

\section{Strong and Weak Cycle Structure Theorems}

Let $A$ be a subset of the positive integers. We consider the class of permutations
that consist entirely of disjoint cycles of lengths in $A$, and denote this by 
$\fancy{P}^{(A,\Z^{\geq 0})}$. Furthermore, if $B \subseteq \Z^{\geq 0}$, we may
consider the subclass $\fancy{P}^{(A,B)} \subseteq \fancy{P}^{(A,\Z^{\geq 0})}$
consisting of only those permutations whose cycle count is found in $B$. That is, any permutation
of cycle count not in $B$, or containing a cycle length not in $A$, are prohibited.  

The following theorems were first proven (presumably) long ago but can be found in
\cite{Flajolet} and also \cite{MarkoRiedel}, and it is commonly noted that the technique
in general was used by Laplace in the late 18th century. The nomenclature is however, ours.

\begin{theorem}\label{strongtheorem}
The Strong Cycle Structure Theorem:

The combinatorial class $\fancy{P}^{(A,B)}$ has associated EGF, $\fancy{P}^{(A,B)}_e(z)
= \beta(\alpha(z))$, where $\beta(z)$ is the EGF associated to $B$ and 
$\displaystyle{\alpha(z) = \sum_{i \in A} \frac{z^i}{i}}$.
\end{theorem}

However, we only need a weaker form in all but one case in this paper:

\begin{theorem}\label{weaktheorem}
The Weak Cycle Structure Theorem:

The combinatorial class $\fancy{P}^{(A,\Z^{\geq 0})}$ has associated EGF, 
$\fancy{P}^{(A,\Z^{\geq 0})}_e(z) = exp(\alpha(z))$, where $\alpha(z)$ is as above:
$\displaystyle{\alpha(z) = \sum_{i \in A} \frac{z^i}{i}}$
\end{theorem}

This is clearly a special case of the Strong Cycle Structure Theorem with $\beta(z) =
1+z + z^2/2!+z^3/3!+z^4/4!+\cdots = e^z$ (the EGF of $\Z^{\geq 0}$). Interestingly, 
if $A = \Z^+$, then 
$\alpha(z) = z + z^2/2 + z^3/3 + z^4/4 + z^5/5 + \cdots = \log\left(\frac{1}{1-z}\right)$, which
provides a verification of the theorem in this special case:
\[
  exp\left(\log\left(\frac{1}{1-z}\right)\right) = \frac{1}{1-z} = 1+z + z^2 + 
  z^3 + z^4 + \cdots,
\]
which is the EGF for the combinatorial class $\fancy{P}$ of all permutations (together
with the unique permutation on $0$ letters), as expected.

Since the proof of the strong version is not fundamentally more difficult than the weak version,
we shall provide a proof of Theorem~\ref{strongtheorem}.  While this has been proven already in~\cite{Flajolet}, 
we feel that a more expository proof is appropriate in this context.  First, a
lemma which proves the case $B = \{k\}$.

\begin{lemma}\label{lem.B_is_k}  
  The combinatorial class $\fancy{P}^{(A,\{k\})}$ has associated EGF, 
  \[
    \fancy{P}_e^{(A, \{k\})}(z) = \frac{1}{k!}\left( \sum_{i \in A} \frac{z^i}{i} \right)^k.
  \]
\end{lemma}
\begin{proof}
  Let $A \subseteq \Z^+$.  For a given cycle-count, $k$, we must only include cycles of lengths found in $A$.
  Begin with an OGF.  If $\pi \in S_n$ has $k$ cycles, then its cycle structure defines a partition
  of $n$ identical objects into $k$ containers, where each container cannot have any number of
  objects that does not occur as a member of $A$.  The OGF that generates this is
  $\left( \sum_{i \in A} z^i \right)^k$, as stated in Section~\ref{statedhere}. 
  Now, we must remember that those objects in the containers
  are {\it not} identical!  Think of each cycle-structure as being a template onto which we
  attach the labels $1, 2, 3, 4, \ldots, n$ in some order.  A priori, this provides a factor
  of $n!$ for each partition of $n$, and so the coefficient of $z^n$ in the above OGF should be
  multiplied by $n!$.  
  The best way to accomplish this is to simply consider our OGF as an EGF:  In
  our OGF, if $C_n$ is the coeffiecient of $z^n$, then as EGF, $n! C_n$ is the
  coefficient of $z^n/n!$.  Now, for each disjoint cycle of length
  $i$, there are $i$ ways of cyclically permuting the labels, each giving rise 
  to an equivalent representaion of the same $i$-cycle.
  Thus, we have over-counted unless we divide each term $z^i$ by $i$.  Finally, each rearrangement
  of the $k$ cycles among themselves gives rise to an equivalent expression for the permuation,
  so we must divide by $k!$, and our EGF for permutations of cycle-count $k$ with cycle-lengths
  in $A$ now has the required form, $\fancy{P}_e^{(A, \{k\})}(z) = \frac{1}{k!}
  \left( \sum_{i \in A} z^i/i \right)^k$.
\end{proof}

The proof of Theorem~\ref{strongtheorem} then follows easily:

\begin{proof}
  Let $A\subseteq \Z^+, B \subseteq \Z^{\geq 0}$.  Categorize all permutations
  in $\fancy{P}$ by cycle-count.
  Only permutations with cycle-counts $k \in B$ will contribute to our total, 
  so by Lemma~\ref{lem.B_is_k},
  \[
    \fancy{P}_e^{(A, B)}(z) = \sum_{ k \in B} \fancy{P}_e^{(A, \{k\})}(z)= \sum_{k \in B}
    \frac{1}{k!}\left( \sum_{i \in A} \frac{z^i}{i} \right)^k
    = \sum_{ k \in B} \frac{\alpha(z)^k}{k!} = \beta(\alpha(z)),
  \]
  since $\sum_{k \in B} z^k/k!$ is the EGF associated to $B$.  The Weak Cycle 
  Structure Theorem then follows as an immediate corollary.  
\end{proof}

\subsection{Probabilities}
In cryptography and other disciplines, we are often concerned with determining
{whe-ther} or not a random permutation has some given property $\phi$. We can calculate
then the OGF of the combinatorial class $\fancy{F}$ of permutations with that property, and
divide term-wise with the same term from the OGF of $\fancy{P}$, the combinatorial
class of all permutations. But this is the same as the coefficients of the EGF of $\fancy{F}$.  

This works for any specific size, but first, it might be difficult to calculate, and second
we might want to know the limit of this probability as the size goes to infinity.

\begin{theorem}\label{hardtheorem}
Let $\fancy{F} \subset \fancy{P}$ be the combinatorial class of permutations with
property $\phi$. Suppose further $\fancy{F}$ has EGF equal to $f(z)$.
Then the limit (as $n$ goes to infinity) of the 
probability that a random permutation of size $n$ has property $\phi$ is given
by 
$$p = \lim_{z\rightarrow 1^-} (1-z)f(z)$$ 
provided that $(1-z)f(z)$ is continuous from the left at $z=1$.
\end{theorem}

\begin{proof}
Let the OGF of $\fancy{F}$ be given by $A_0 + A_1 z + A_2 z^2 +  A_3 z^3 +  
A_4 z^4 +  A_5 z^5 + \cdots$.  Consider the following function
$$
g_n(z) = \frac{A_0}{0!} + \sum_{1\leq i\leq n} \left ( \frac{A_i}{i!}-
\frac{A_{i-1}}{(i-1)!} \right ) z^i,
$$
which when evaluated at $z=1$, the sum telescopes,
$$
= \frac{A_0}{0!} + \left( \frac{A_1}{1!}-\frac{A_{0}}{0!} \right ) (1)
+ \left( \frac{A_2}{2!}-\frac{A_{1}}{1!} \right ) (1)^2 + \cdots
+ \left( \frac{A_n}{n!}-\frac{A_{n-1}}{(n-1)!} \right ) (1)^n = \frac{A_n}{n!}.
$$
%
Thus $g_n(1)$ is the desired probability, for size $n$.

The limit $g(z) = \lim_{n\rightarrow \infty} g_n(z) = 
\frac{A_0}{0!} + \sum_{i\ge 1} \left ( \frac{A_i}{i!}-\frac{A_{i-1}}{(i-1)!} \right ) z^i$
does not necessarily exist for all $z$, but when it does, we have
\begin{eqnarray*}
g(z) = \lim_{n\rightarrow \infty} g_n(z)& = & \lim_{n\rightarrow \infty} \frac{A_0}{0!} + 
\left ( \sum_{i = 1}^n  \frac{A_i}{i!}z^i \right )
- \left ( \sum_{i = 1}^n \frac{A_{i-1}}{(i-1)!}z^i \right ) \\
\ignore{&=&  \lim_{n\rightarrow \infty} \frac{A_0}{0!} + 
\left ( \sum_{i = 1}^n  \frac{A_i}{i!}z^i \right )
- \left ( \sum_{j = 0}^n \frac{A_j}{j!}z^{j+1} \right ) \\}
& = &  \lim_{n\rightarrow \infty} \left ( \sum_{i = 0}^n \frac{A_i}{i!}z^{i} \right )
- z \left ( \sum_{j = 0}^n \frac{A_j}{j!}z^{j} \right )\\
& = & (1-z) \lim_{n\rightarrow \infty} \left ( \sum_{i = 0}^n \frac{A_i}{i!}z^{i}
\right ) = (1-z)f(z)
\end{eqnarray*}

Thus $p=\lim_{n\rightarrow \infty} g_n(1) = \lim_{n\rightarrow \infty} 
\lim_{z\rightarrow 1^-} g_n(z) = \lim_{z\rightarrow 1^-} (1-z)f(z)$.

Note, we implicitly assumed that $g(z)$ is continuous (from the left)
near $z=1$ in order to reverse the order of the limits in the last step, but this will be
the case in all of our examples.
\end{proof}

\subsection{Expected Values}

While OGFs and EGFs are very useful for the study of a one-parameter family of constants,
$A_0, A_1, A_2, A_3, \ldots$, we often wish to work with a two-parameter family, $\{A_{s,t} 
\}_{s,t \geq 0}$.  This is accomplished using {\it double} generating functions.  The double
OGF, $A(y,z)$ of a two-parameter family of constants, $\{A_{s,t}\}$ is defined to be the 
formal sum: 
\[
  A(y,z) = \sum_{s = 0}^{\infty} \sum_{t=0}^{\infty} A_{s,t}y^sz^t,
\]
and the EGF $A_e(y,z)$ is defined to be the formal sum:
\[
  A_e(y,z) = \sum_{s = 0}^{\infty} \sum_{t=0}^{\infty} \frac{A_{s,t}}{(s+t)!}y^sz^t.
\]

For our purposes, we will be interested in a combinatorial class of permutations categorized
not only by the order of the symmetric group $S_n$ in which the permutation lies, but also by the
number of fixed points that the permutation possesses.

\begin{theorem}\label{thm.exp_value}
Let $\fancy{F} \subset \fancy{P}$ be a combinatorial class of permutations
with double EGF $a(y,z)$, where the coefficient of $y^s z^t/(s+t)!$
is the number of permutations $\pi$ with property $\phi_s$ such that $\pi \in S_{s+t}$.  Then
the limit (as $n = s+t$ goes to infinity) of the expected value of $s$ such that a random
permutation of size $n$ satisfies $\phi_s$ is given by:
\[
  \lim_{z \to 1^-} (1-z) a_y(z,z)
\]
provided $(1-z)a_y(z,z)$ is convergent and continuous from the left at $z=1$. 
\end{theorem}

\begin{proof}
  Let $a(y,z) = \sum_{s \geq 0}\sum_{t \geq 0} y^s z^tA_{s,t}/(s+t)! $.  The coefficient of
  $y^s z^t$ is the probability that a random permutation of $S_{s+t}$ has property $\phi_s$,
  by construction.
  Consider the partial derivative with respect to $y$:  
  \[
    a_y(y,z) = \sum_{s \geq 0}\sum_{t \geq 0} \frac{sA_{s,t}}{(s+t)!} y^{s-1}z^t. 
  \]
  The probabilities are now multiplied by the corresponding value of $s$.  Now, letting
  $y = z$ produces:
  \[
    a_y(z,z) = \sum_{s \geq 0}\sum_{t \geq 0} \frac{sA_{s,t}}{(s+t)!} z^{s+t-1}
      = \sum_{n \geq 0} \left(\sum_{s+t = n} \frac{sA_{s,t}}{n!}\right) z^{n-1} .
  \]
  Thus, $a_y(z,z)$ is the OGF that computes the expected value of $s$ such that a random
  permutation of size $n$ satisfies $\phi_s$ (shifted by one degree).  Using the same technique 
  as in the proof of Thm~\ref{hardtheorem}, we find that
  \[
    \lim_{z \to 1^-}   (1-z) a_y(z,z) =  \lim_{n \to \infty}
    \left(\sum_{s+t = n} \frac{sA_{s,t}}{n!}\right).
  \]
\end{proof}

\section{Corollaries}
Theorem~\ref{hardtheorem} is exploited extensively in a paper by Marko R. Riedel
dedicated to random permutation statistics, but in a different context (see  \cite{MarkoRiedel}).

\begin{corollary}\label{CoroC3}
The probability that a random permutation (in  the limit as the size grows
to infinity) does not contain cycles
of length $k$ is given by $e^{-1/k}$.
\end{corollary}

\begin{proof}
The set $A$ of allowable cycle lengths is $\integer^+ - \set{k}$,
and so has EGF given by artificially removing the term for $k$ from the EGF of $\fancy{O}$:
\ignore{$$ z + z^2 + 2!z^3 + \cdots + (k-2)!z^{k-1} + 0 + k!z^{k+1} + (k+1)!z^{k+2} 
+ (k+2)!z^{k+3} + \cdots$$
and therefore EGF}
$$ z + \frac{z^2}{2} + \frac{z^3}{3} + \cdots +  
\frac{z^{k-1}}{k-1} + 0 + \frac{z^{k+1}}{k+1} + \frac{z^{k+2}}{k+2} 
+ \cdots = \log \left ( \frac{1}{1-z} \right ) - \frac{z^k}{k},$$
and thus by the Weak Cycle Structure Theorem, the combinatorial class
in question has EGF equal to 
$$ a(z) = exp \left (  \log \left ( \frac{1}{1-z} \right ) - \frac{z^k}{k}
\right ) = \frac{1}{1-z} e^{-z^k/k}$$

Thus the probability of a random permutation (as the size tends toward infinity)
not having any cycles of length $k$ is given by
$\lim_{z\rightarrow 1^-} (1-z)a(z) = e^{-1/k}$
\end{proof}

\paragraph{Note: On the Precision of these estimations:} 
This result means that $p\to e^{-\frac{1}{k}}$ when $N\to \infty$.
What about when $N=2^{32}$? We can answer this question easily by
observing that the Taylor expansion of the function $a(z)$ is the EGF and
therefore gives all the \emph{exact} values of $A_n/n !$.
For example when $k=4$ we computed the Taylor expansion of $g(z)$ at order 201,
where each coefficient is a computed as a ratio of two large integers.
This takes less than a second with the computer algebra software \emph{Maple} \cite{Maple}.
The results are surprisingly precise: the difference between
the $A_{200}/200!$ and the limit is less than $2^{-321}$.
Thus convergence is very fast and even for very small permutations (on 200 elements).

Returning to the proving of corollaries,
let us define $\fancy{P}^{\overline{A}}=\fancy{P}^{(\integer^+-A,
\integer^{\geq 0})}$ and find its EGF.

\begin{lemma}\label{prohibcycles}
The EGF of $\fancy{P}^{\overline{A}}$ is given by $exp\left ( f(z) \right )$,
where 
$$f(z) = \sum_{i \not\in A} z^i/i = \log \left ( \frac{1}{1-z} \right ) -  
\sum_{i \in A} z^i/i $$
\end{lemma}

\begin{proof}
Because  $\fancy{P}^{\overline{A}}=\fancy{P}^{(\integer^+-A,\integer^{\geq 0})}$
we can use the Weak Cycle Structure Theorem.
The EGF of the combinatorial class of cycles with size from the
set $\integer^+-A$ is given by that of $\fancy{O}$ (the class of all cycles)
with the ``forbidden lengths'' artificially set to zero, namely
$$ \sum_{i \in (\integer^+ - A)} z^i/i = \sum_{0<i \not\in A} z^i/i = 
\log \left ( \frac{1}{1-z} \right ) -  \sum_{i \in A} z^i/i$$
The correct answer follows.
\end{proof}

\begin{corollary}\label{PropC5}
Let $A$ be a subset of the positive integers.
The probability that a random permutation (in  the limit as the size grows
to infinity) does not contain cycles of length in $A$ is:
$$\prod_{i\in A} e^{-1/i} = e^{-\sum_{i\in A} 1/i} $$
\end{corollary}

\begin{proof}
Using Lemma~\ref{prohibcycles} we obtain an EGF of
$$ exp\left ( \log \left ( \frac{1}{1-z} \right ) - \sum_{i\in A} z^i/i \right)
= \frac{1}{1-z} \prod_{i\in A} e^{-z^i/i}$$
then multiplying by $(1-z)$ and taking the limit as $z\rightarrow 1$ gives the desired result.
\end{proof}

This offers confirmation of Corollary~\ref{CoroC3} when substituting $A=\set{k}$.
A permutation with no fixed points is called a derangement. Using a similar
strategy, we can calculate the probability of a derangement.

\begin{corollary}\label{derangedcorollary}
Let $\pi$ be a permutation taken at random from $S_n$. The
probability that $\pi$ is a derangement is $1/e$ in the
limit as $n \rightarrow \infty$.
\end{corollary}

\begin{proof}
Just apply Corollary~\ref{PropC5} to the case of cycle length 1. 
\end{proof}

Suppose we wish to consider if a permutation has exactly $t$ cycles of length
from a set $C \subset \integer^+$, 
in other words, all the other cycles are of length not found in $C$. 
In that case, we can consider such a permutation $\pi$ as a product of $\pi_A$
and $\pi_B$ such that $\pi_A$ has only $t$ cycles of length found in $A$, and
nothing else, and $\pi_B$ has only cycles of length not found in $A$. This is
termed by Flajolet and Sedgewick as a ``labelled product''\footnote{A labelled product
can be thought of as follows. If the EGF of $a(z)=b(z)c(z)$, where $b$ and $c$ are
also EGFs, then $a(z)=\sum_{k=0}^{k=n} {n \choose k}b_kc_{n-k}$. Here, after 
building our combinatorial object in class $a$ of size $n$ out of `an object' from $b$ of size $k$,
and `an object' from $c$ of size $n-k$, we must then attach $k$ of the $n$ labels to the
former, and attach the remaining $n-k$ labels to the latter. There are precisely ${n \choose k}={n \choose {n-k}}$
ways to do that.} and and a discussion
of that is found in Section II.2 in \cite{Flajolet}. The EGF of a labelled
product is merely the product of the EGFs.

\begin{theorem}\label{thm.t-fixed-pts}
  Let $\pi$ be a permutation taken at random from $S_n$.  The probability that $\pi$ has $c$
  fixed points is $\frac{1}{c! e}$.
\end{theorem}
\begin{proof}
  Consider $\pi = \pi_A\pi_B$, where $\pi_A$ consists of exactly $c$ fixed points, and $\pi_B$
  is a derangement of the remaining $n - c$ points.  We must compute the labelled product
  $f(z) = \fancy{P}_e^{(\{1\}, \{c\})} \cdot \fancy{P}_e^{( \Z^+ - \{1\}, \Z^{\geq 0})}$.  Thus,
  by the Strong and Weak Cycle Structure Theorems,
  \[
    f(z) = \frac{z^c}{c!} exp\left( \log\left(\frac{1}{1-z}\right) - z \right) = 
    \frac{z^c}{(1-z)c!} e^{-z}
  \]
  An application of Thm~\ref{hardtheorem} provides the result:
  \[
    \lim_{z \to 1^-} (1-z)f(z) = \lim_{z \to 1^-}\frac{z^c}{c!} e^{-1} = \frac{1}{c! e}
  \]
\end{proof}

\subsection{On Cycles in Iterated Permutations}
\label{RPStats3}

\begin{theorem}\label{fixedfromdivisor}
Let $\pi$ be a permutation in $S_n$.
A point $x$ is a fixed point for $\pi^k$ if and only if
$x$ is a member of a cycle of length $i$ in $\pi$, for some
positive integer $i$ dividing $k$.
\end{theorem}

\begin{proof}
Write $\pi$ in disjoint cycle notation, and then $x$ appears
in only one cycle (hence the name ``disjoint.'') Call this
cycle $\psi$. Since all other cycles do not contain $x$, 
then $\pi^m(x)=\psi^m(x)$ for all integers $m$. Of course,
$\psi$ is of order $i$ in $S_n$, thus $\psi^i=id$, the 
identity element of $S_n$.

If $x$ is in a cycle of length $i$ then that means
that $i$ is the smallest positive integer such that
$\psi^{i}(x)=x$. Write $k=qi+r$ with $0\le r < i$. Then
$$x= \psi^k(x) = \psi^r( \psi^{iq}(x) ) = \psi^r( (\psi^i)^q(x) ) = 
\psi^r( id^q(x) ) = \psi^r( id(x) ) = \psi^r(x) $$
so $\psi^r(x)=x$
but we said that $i$ is the least positive integer such that
$\psi^i(x)=x$ and $r<i$. The only way this is possible is if
$r$ is not positive, i.e. it is zero. Thus $k=qi$ or $i$ divides
$k$.

There reverse assumes that $i$ divides $k$ so write $iq=k$ then
$$ \psi^k(x) = \psi^{iq}(x) = (\psi^i)^q(x) = (id)^q(x) = id(x) = x$$ 
\end{proof}

\paragraph{An Example}
Before we continue, observe what happens to a cycle of $\pi$
when evaluating $\pi^2$. First, if the cycle is of odd length,
$$ (x_1, x_2, \ldots, x_{2c+1}) \mapsto
(x_1, x_3, x_5, \ldots, x_{2c+1}, x_2, x_4, x_6, \ldots, x_{2c})$$
but if the cycle is of even length,
$$ (x_1, x_2, \ldots, x_{2c}) \mapsto
(x_1, x_3, x_5, \ldots, x_{2c-1})(x_2, x_4, x_6, x_8, \ldots, x_{2c})$$


One can rephrase Theorem~\ref{fixedfromdivisor} as follows:

\begin{corollary}\label{divisorcorollary}
Let $\pi$ be a permutation from $S_n$.
Let $k$ be a positive integer, and let the set of positive
integer divisors of $k$ be $D$. Then the set of fixed points
of $\pi^k$ is precisely 
the set of points under $\pi$ in cycles of length found in $D$.
\end{corollary}

\subsection{Limited Cycle Counts}\label{limcyclecount}

\ignore{
The EGF of the combinatorial class of permutations with cycles only from $A$
is $\fancy{P}^{(A, \integer^+)}$ as before. However, we must then restrict ourselves
to those having \emph{exactly} $t$ such cycles.

\begin{theorem}\label{specialtheorem}
Let $\pi$ be a permutation taken at random from $S_n$.
The probability that a random permutation (in the
limit as $n\rightarrow \infty$) 
has precisely $t$ cycles of length found in $A\subset \integer^+$
is given by
$$\frac{1}{t!} e^{-\sum_{i \in A} 1/i}\left ( \sum_{i\in A} 1/i \right )^t$$
\end{theorem}

\begin{proof}
Consider $\pi = \pi_A\pi_B$ as in the discussion above, where $\pi_A$
consists of exactly $t$ disjoint cycles with cycle types found only in $A$, and
$\pi_B$ consists of any number of disjoint cycles of types not found in $A$.
We must compute the labelled product $\fancy{P}_e^{(A, \{t\})}\cdot
\fancy{P}_e^{(\overline{A}, \Z^+)}$
By Lemma~\ref{lem.B_is_k}, $\fancy{P}_e^{(A, \{t\})}(z) = \frac{1}{t!}
\left( \sum_{i \in A}\frac{z^i}{i} \right)^t$, and by Lemma~\ref{prohibcycles},
$\fancy{P}_e^{(\overline{A}, \Z^+)}(z) = \frac{1}{1-z}\prod_{i \in A}
e^{-z^i/i}$.  Thus, our EGF is:
\[
  a(z) = \frac{1}{t!}\left( \sum_{i \in A}\frac{z^i}{i} \right)^t\left(
  \frac{1}{1-z}e^{-\sum_{i \in A}z^i/i}\right)
  = \frac{1}{(1-z)t!}e^{-\sum_{i \in A}z^i/i}\left( \sum_{i \in A}\frac{z^i}{i}
  \right)^t
\]
Then we see (by invoking Theorem~\ref{hardtheorem})
that $(1-z)a(z)$ gives the desired probability as $z \rightarrow 1$:
\[
  \lim_{z \to 1^-} (1-z) \cdot \frac{1}{t!}\left( \sum_{i \in A}\frac{z^i}{i} 
  \right)^t\left(\frac{1}{1-z}e^{-\sum_{i \in A}z^i/i}\right)
  = \lim_{z \to 1^-}\frac{1}{t!}\left( \sum_{i \in A}\frac{z^i}{i} \right)^t
  e^{-\sum_{i \in A}z^i/i}
\]
\end{proof}

\begin{corollary}\label{coolcorollary}

[FIXME!!!]

Let $\pi$ be a permutation taken at random from $S_n$.
The probability that $\pi^k$ has exactly $t$ fixed points (taken in the
limit as $n$ goes to infinity) is given by
$$\frac{1}{t!} e^{-\sigma(k)/k}\left ( \frac{\sigma(k)}{k} \right )^t$$
\end{corollary}

\begin{proof}
A point $x$ is a fixed point under $\pi^k$ if and only if $x$ is a member of a
cycle of order dividing $k$ under $\pi$, via Corollary~\ref{fixedfromdivisor}. 
Thus if $\pi$ has exactly $t$ cycles of length dividing $k$, then $\pi^k$ has 
exactly $t$ fixed points.  Thus, by Thm.~\ref{specialtheorem}, with $A$ being
the set of all divisors of $k$, the probability works out to be:
\[
  \frac{1}{t!} e^{-\sum_{i \in A} 1/i}\left ( \sum_{i\in A} 1/i \right )^t
  = \frac{1}{t!} e^{-\sum_{i | k} 1/i}\left ( \sum_{i |k} 1/i \right )^t
  = \frac{1}{t!} e^{-\sigma(k)/k}\left ( \frac{\sigma(k)}{k} \right )^t.
\]
Note, the last equality follows from Lemma~\ref{sigmalemma}.
\end{proof}

Furthermore, note that substituting $t=0$ into the above
yields the same result as Corollary~\ref{derangedcorollary}, as expected.

\begin{theorem}
Let $k$ be a positive integer, and $\pi$ a permutation from $S_n$.
The expected number of fixed points of $\pi^k$ is $\sigma(k)/k$, taken in the
limit as $n \rightarrow \infty$.
\end{theorem}

\begin{proof}
[FIXME!!!]
\end{proof}
}

\begin{theorem}\label{tautheorem}
  Let $k$ be a positive integer, and $\pi$ a permutation from $S_n$.
  The expected number of fixed points of $\pi^k$ is $\tau(k)$, taken in the
  limit as $n \rightarrow \infty$. Note, $\tau(k)$ is the number of positive
  integers dividing $k$.
\end{theorem}
\begin{proof}
  We shall construct a double EGF, $a(y,z)$, where the coefficient of $y^s z^t$ is the probability
  that the $k^{th}$ power of a random permutation of $S_{s+t}$ has $s$ fixed points.  Let $\pi$ be 
  a permutation taken at random from $S_n$.
  A point $x$ is a fixed point under $\pi^k$ if and only if $x$ is a member of a
  cycle of order dividing $k$ under $\pi$, via Corollary~\ref{fixedfromdivisor}. 
  Thus $\pi^k$ has exactly $t$ fixed points if and only if $\pi = \pi_A\pi_B$, where
  $\pi_A \in S_t$ consists only of cycles of length dividing $k$, and $\pi_B \in
  S_{n-t}$ consists only of cycles of length {\it not} dividing $k$.  Let $D_k$ be
  the set of all positive divisors of $k$.  The double
  EGF that counts the number of such permutations $\pi_A\pi_B$ will be given by the
  labelled product $\fancy{P}_e^{(D_k, \Z^{\geq 0})}(y) \cdot \fancy{P}_e^{(\Z^+ - D_k, 
  \Z^{\geq 0})}(z)$. By the Weak Cycle Structure Theorem and Lemma~\ref{prohibcycles},
  we obtain:
  \begin{eqnarray*}
    a(y,z) &=& exp\left( \sum_{i | k} \frac{y^i}{i} \right) exp\left( \log\left(\frac{1}{1-z}\right)
    - \sum_{i | k} \frac{z^i}{i} \right)\\
    &=& exp\left(\log\left(\frac{1}{1-z}\right)\right)exp\left(
    \sum_{i | k} \frac{y^i}{i} - \sum_{i | k} \frac{z^i}{i} \right)\\
    &=& \frac{1}{1-z} exp\left( \sum_{i | k} \frac{y^i-z^i}{i} \right).
  \end{eqnarray*}
  Theorem~\ref{thm.exp_value} provides the correct expected value.  First observe that
  \[
    a_y(y,z) = \frac{1}{1-z} exp\left( \sum_{i | k} \frac{y^i-z^i}{i} \right)\sum_{i | k} y^{i-1}.
  \]
  Then $a_y(z,z) = \frac{1}{1-z} exp(0)\sum_{i | k} z^{i-1}$.  Finally,
  \[
    \lim_{z \to 1^-} (1-z) a_y(z,z) = \lim_{z \to 1^-} \sum_{i | k} z^{i-1} = 
    \sum_{i | k} 1 = \tau(k).
  \]
\end{proof}

\section{Application to Keeloq}

\subsection{What is Keeloq?}

Keeloq is a block cipher, with 32-bit plaintext and ciphertext blocks and a 64-bit key.
It has been used in the remote keyless entry systems of many manufacturers of
automobiles, and several papers have been written about it 
\cite{KeeLoqDunkelman}
\cite{BogdanovKeeLoq}
\cite{BogdanovKeeLoq2}
\cite{BogdanovKeeLoq3}
\cite{KeeLoq0}
\cite{KeeLoq}
\cite[Ch. 2]{BardDis}.
\cite{KeeLoqTatra}. It has 528 rounds,
which is unusually high, and this can be written $528=8\times 64 + 16$, a decomposition
whose utility will be apparent shortly. Each 
round is like a stream cipher, in the sense that the internal state is a 32-bit register,
and is shifted
one bit, and a new bit is introduced. The new bit is a function of certain bits of the 
internal state, and a single bit of the key, via a map described by a
cubic polynomial over GF(2), see for example \cite{monograph} \cite[Ch. 2]{BardDis}. 
The initial
value of the internal state is the plaintext, and the final value is the ciphertext.
For completeness, the cipher specification is given in Figure 1.

\begin{figure}
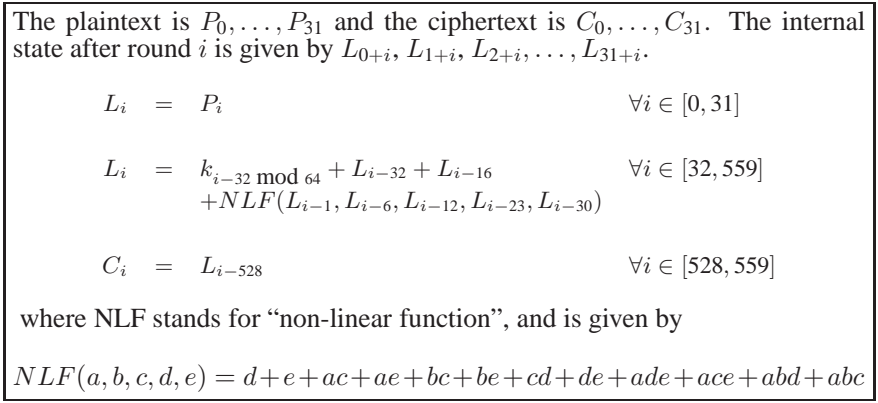

\begin{center}
\framebox{\begin{minipage}{0.9\textwidth}
The plaintext is $P_0, \ldots, P_{31}$ and the ciphertext is $C_0, \ldots, C_{31}$.
The internal state after round $i$ is given 
by $L_{0+i}$, $L_{1+i}$, $L_{2+i}, \ldots, L_{31+i}$.
{\small$$\begin{array}{rcll}
L_i & = & P_i & \forall i \in [0, 31]\\
\\
L_i & = & k_{i-32 \mbox{ mod }64} +  L_{i-32} +  L_{i-16}  & \forall i \in [32, 559]\\
& & +  
NLF(L_{i-1}, L_{i-6}, L_{i-12}, L_{i-23}, L_{i-30})\\
\\
C_i & = & L_{i-528} & \forall i \in [528, 559]
\end{array}$$}
where NLF stands for ``non-linear function'', and is given by
$$
NLF(a,b,c,d,e)=d +  e +  ac +  ae +  bc +  be
+  cd +  de +  ade +  ace +  abd +  abc
$$
\end{minipage}}
\end{center}
\caption{The Specification of Keeloq}
\end{figure}

Also, because each round only uses 1 bit of the key (and they are used in sequence), 
then after 64 rounds, the entire key has been used. 
Therefore, it makes sense to define $f_k$, a function which represents those
64 rounds. Each additional 64 rounds behaves identically.
It turns out that $f_k$ is a permutation. The remaining 16 rounds
are written as $g_k$, which is also a permutation. Of course if either $f_k$ or
$g_k$ were not permutations, then the block cipher would not be uniquely decodable.

Thus we can write
$ g_k(f^{(8)}_k(p))=E_k(p)$
and this motivated the authors' initial interest in iterated permutations. 
Also it is noteworthy that only 16 bits of the key are used by $g_k$, thus
only 16 bits of the key need be known or guessed to use $g^{-1}_k$ to 
``peel off'' or ``undo'' these 16 rounds, leaving us with $f^{(8)}_k$, the
eighth iterate of a permutation.

\subsection{Bard's Dissertation Attack}

This attack assumes some portion of the code-book is available. So long as two
fixed points are found, the attack can succeed. One can show that if there are
two plaintexts that are fixed on the first 64 rounds of the encryption, i.e.
$f(p_1)=p_1$ and $f(p_2)=p_2$, then this is sufficient information to perform
an algebraic cryptanalysis, see \cite[Ch. 2]{BardDis} \cite[Ch. 3]{monograph}. 
One writes polynomials for those two
equalities and uses SAT-solvers to solve them, see \cite{SATconv} \cite[Ch. 6]{monograph}.

The question becomes how to obtain those pairs.
First, the part of the key used in $g_k$, which is 16 bits in length, is simply guessed.
This has success probability $2^{-16}$. Then $g^{-1}_k$ can be used. This
allows for $(p,c)$, the plaintext-ciphertext pairs in the codebook to be replaced
by $(p,g^{-1}_k(c))$ which are now actually $(p,f^{(8)}_k(p))$. These are points
fixed by $f^{(8)}_k$ and so by Corollary~\ref{fixedfromdivisor}, 
they are points of order $\set{1,2,4,8}$
for $f_k$. Thus, the fixed points of $f_k$, which are useable for the cryptanalysis,
are a subset of those for $f^{(8)}_k$, which we can find.

\begin{theorem}\label{magictheorem}
Let $\pi$ be a random permutation from $S_n$. The probability that
$\pi$ has $c_1$ fixed points and $c_2$ cycles of lengths 2, 4, or 8,  
is given by
$$ \frac{1}{c_1!c_2!}\left ( \frac{7}{8} \right )^{c_2} e^{-15/8} $$
\end{theorem}

\begin{proof}
Note that the set of permutations
on $n$ elements, with $c_1$ fixed points, and $c_2$ cycles of length
2, 4, or 8, can be thought of as a triple labelled product. The first item in the
product is from $\fancy{P}^{(\set{1},c_1)}$, the second item from $\fancy{P}^{(\set{2,4,8}, c_2)}$,
and the third item from $\fancy{P}^{\overline{\set{1,2,4,8}}}$. We must now calculate
the EGF.

The first item has $\alpha(z)=z$, and $\beta(z)=z^{c_1}/{c_1!}$, for an EGF
of $\beta(\alpha(z))=z^{c_1}/{c_1!}$. The second item has
$\alpha(z)=z^2/2+z^4/4+z^8/8$, and $\beta(z)=z^{c_2}/{c_2!}$, therefore
an EGF of 
$\beta(\alpha(z))=\frac{1}{c_2!}\left [ z^2/2 + z^4/4 + z^8/8 \right ]^{c_2}$.
Finally, the third item has EGF given by Lemma~\ref{prohibcycles},
$$ exp \left ( \log \left ( \frac{1}{1-z} - \sum_{i | 8} z^i/i \right ) \right )
= \frac{1}{1-z} exp \left (  -\sum_{i | 8} z^i/i \right ) $$
giving a final, total EGF of
$$ \frac{z^{c_1}}{(1-z)c_1!c_2!} \left [ \frac{z^2}{2} + \frac{z^4}{4} + 
\frac{z^8}{8} \right ]^{c_2} exp \left ( - \sum_{i|8} z^i/i \right )$$
Multiplying by $1-z$ and taking the limit as $z\rightarrow 1^-$, via Theorem~\ref{hardtheorem} we obtain
$$ \frac{1}{c_1!c_2!} \left [ \frac{1}{2} + \frac{1}{4} + \frac{1}{8} \right ]^{c_2} exp ( -\sigma(8)/8 ) = \frac{1}{c_1!c_2!}(7/8)^{c_2}e^{-15/8}$$ 
\end{proof}

The method requires $c_1 \ge 2$, otherwise the attack fails. This can be easily
calculated as $1-\Pr \{ c_1=0 \} -\Pr \{ c_1=1 \} \approx 0.2642$ probability of success.

Second, suppose that $\eta$ is the fraction of the code-book available. Then any given fixed point is found
with probability $\eta$ in the known part of the code-book, and so at least two will be found with probability
$$ 1 - {c_1 \choose 0}\eta^0(1-\eta)^{c_1} - {c_1 \choose 1}\eta^1(1-\eta)^{c_1-1} = 
1 - (1-\eta)^{c_1-1}\left [ 1 - (c_1+1)\eta\right ]$$
and so the following $\eta$ and success probabilities can be found,
generated by Theorem~\ref{thm.t-fixed-pts} and
listed in Table~\ref{bardsuccess}. Note, these are absolute probabilities,
not probabilities given $c_1\ge 2$.

\begin{table}
{\begin{tabular}{ccccccccccc}
$\eta$ & 10\% & 20\% & 30\% & 40\% & 50\% \\
Success & 0.47\% & 1.75\% & 3.69\% & 6.16\% & 9.02\% \\
\\
$\eta$ & 60\% & 70\% & 80\% & 90\% & 100\% \\
Success & 12.19\% & 15.58\% & 19.12\% & 22.75\% & 26.42\%
\end{tabular}}
\label{bardsuccess}
\caption{Success Probabilities of Bard's Dissertation Attack}
\end{table}

Using Maple, one can also calculate exactly when the probability of having the two fixed points
in the $\eta$ fraction of the code-book is one-half. This is at $\eta=63.2\%$ remarkably
close to the empirical calculation in \cite[Ch. 2]{BardDis}.

Note that while finding two fixed points of $f_k$ is enough to break the cipher, using SAT-solvers as
noted above, the fixed points of $f^{(8)}_k$ are still an annoyance. Our post-processed code-book will 
have all the fixed points of $f^{(8)}_k$ in it, and at worst we must try all pairs. 

If $\pi$ has $c_1$ fixed points, and $c_2$ cycles of length 2, 4, or 8, then $\pi^8$ has
at most $c_1+8c_2$ fixed points, as each cycle of length 2 produces 2, of length 4
produces 4, and of length 8 produces 8. Thus of the $c_2$ cycles of length 2, or 4, or 8,
at most $8c_2$ fixed points are produced. This means in the code-book we have
at most $c_1+8c_2$ fixed points, or $(c_1+8c_2)(c_1+8c_2-1)/2$ pairs of them.
At absolute worst, we have to check all of them. The expected value of the number
of pairs, given $c_1 \ge 2$ can be calculated with Maple, and is
$113/2-105/e\approx 17.87$. As each pair takes less than a minute, this is
not the rate-determining step.
\ignore{
Using the above
theorem, when give that $c_1>1$, the expected value of $c_2\approx0.8750$ and 
$c_1\approx=0.6321$,
so this will not be a burden. }

The post-processing of the code-book will take much more 
time, $\eta2^{32}$ Keeloq encryptions, but this is still much 
smaller than brute-forcing the $2^{64}$ keys.
 
\subsection{The Courtois-Bard-Wagner Attack}
Again, in this attack (first published in \cite{KeeLoq}), 
we will iterate over some portion of the code-book. One property of the
cipher Keeloq, is that only one bit is changed per round. Thus the last sixteen rounds, represented
by $g_k(x)$, only affect sixteen bits of the ciphertext. Thus, if $x$ is a fixed point of $f^{(8)}_k$, then
48 out of the 64 bits will match, compared between the plaintext and the ciphertext. One can easily
scan for this property.

This matching property will always occur for a fixed point of $f^{(8)}_k$, but it happens by coincidence
with probability $2^{-16}$. Therefore, the number of code-book entries with this property will be
the number of fixed points of $f^{(8)}_k$, plus an expected $2^{-16}2^{32}=2^{16}$ ``red herrings''. 
What is remarkable, is that \cite{KeeLoq} contains a formula for the 16 key bits that would cause
the effect if it were because the plaintext were a fixed point (i.e. not a coincidence). Therefore,
each code-book entry with the matching property can be tagged with a 16-bit potential sub-key.

As it turns out, the 16-sub key, as well as any single plaintext-ciphertext pair that is a fixed point
of $f_k$, not merely of $f^{(8)}_k$, is enough to mount an algebraic attack. Thus we have the following
steps. Let $c_3$ denote the number of fixed points of $f^{(8)}$.
\begin{enumerate}
\item Check all $2^{32}$ code-book entries for the matching property.
\item Of these (roughly $2^{16}+c_3$) plaintext-ciphertext pairs, 
compute the sub-key that they imply.
\item For each plaintext-ciphertext pair with the property, set up an algebraic cryptanalysis problem
with the one pair, assuming it is a fixed point of $f$, and assuming the sub-key is correct.
\item If an answer is obtained, verify assumptions. If assumptions turned out to be false, or
if the problem is ``unsatisfiable'', go to Step~3.
\end{enumerate}

Sorting upon this sub-key between Step~2 and Step~3 
would reveal which are the likely pairs, as the same sub-key will
tag all the fixed points of $f_k$ and $f^{(8)}_k$. We expect each of the $2^{16}$ ``red-herrings''
to be tagged with
uniformly randomly distributed potential sub-keys. Therefore, in the first very few Step~3 and
Step~4 executions, we would obtain the key.

What is needed for success? First, that $f_k$ have at least one genuine fixed point. This occurs
with probability $1-1/e$, as proven in Corollary~\ref{derangedcorollary}, and is roughly $0.6321$.
Second, the expected amount of work in Step~1 is at most $2^{32}$ Keeloq Encryptions, and
a more precise estimate is found in \cite{KeeLoq}.
Third, Step~2 is negligible. Fourth, for Step~3 and Step~4, we must execute these stages
for each potential sub-key. Given the model of the previous attack,
and using Theorem~\ref{magictheorem}, we can obtain a bound on 
the expected number of repetitions of
Steps~3 and 4. This is upper-bounded by the expected value of $c_1+8c_2$ given that $c_1>0$.
Using Maple, this comes to $113/2-46/e\approx 39.58$, the difference being that we now allow $c_1=1$, which
was previously forbidden. Of course, without the sorting explained in the previous
paragraph, the expected number of Step~3 and Step~4 executions would be around $2^{15}$.

\section{Highly Iterated Ciphers}\label{highlyiterated}
Here we present two attacks, which while no where near practical feasibility, present
surprising results that the authors did not anticipate.

Suppose there were three na\"ive cryptography students, who choose to use
3-DES iterated\footnote{Since the brute
force attack is the optimal attack known at this time, it is perhaps not completely 
unreasonable. The classic UNIX implementations encrypt with a variant of DES 25 times, for
example \cite[Ch. 8]{coolbook}.} approximately one million times, because they are told that this
will slow down a brute force attacker by a factor of one million.
Alice will choose 1,000,000 iterations, Bob will choose
1,081,079 iterations and Charlie will choose 1,081,080 iterations. Intuitively, one
would not expect these three choices to have significantly different security consequences.

However, assuming that the 3-DES cipher for a random key behaves like a randomly
chosen permutation from $S_{2^{64}}$, these permutations will have
$$ \tau(1,000,000)=49\hspace{0.5in} \tau(1,081,079)=2\hspace{0.5in}\tau(1,081,080)=256 $$
fixed points which allows for the following distinguisher attack. It is noteworthy that Charlie's
number is \ignore{the product of the first 8 prime numbers, and thus}the lowest positive integer $x$ to
have $\tau(x)=256$, while Bob's number (only one less) is prime, and thus has $\tau(x-1)=2$.
This enables the dramatic difference in vulnerability to the attack.

In a distinguishing attack, the attacker is presented either with a cipher, or with 
a random permutation from the set of those with the correct domain.
Randomly iterate through $1/64$ of the plain-space. If a fixed point is found, 
guess that one is being given a user cipher. If no fixed point is found, guess random.

In the case of Alice's implementation,
there will be an expected value of $\approx 0.766$ fixed points. In the case of Bob's, 
$1/32$ expected fixed points. In the
case of Charlie's, $4$ expected fixed
points. A random permutation would have $1/64$ expected fixed points. Thus, we can see
that Charlie's would be easily distinguishable from a random permutation, but Bob's much
less so. Against Alice, the attack could definitely still be mounted but with an intermediate
probability of success. 
To make this notion precise, we require the probability distribution of the number of 
fixed points of $\pi^k$. In fact, one can prove the following

\begin{theorem}
Let $\pi \in S_n$ be a permutation chosen at random, then
the $c^{th}$ term of the following EGF
$$exp \left ( \sum_{i|k} \frac{y^i-1}{i} \right )$$
is the probability that $\pi^k$ has exactly $c$ fixed points.
\end{theorem}
\begin{proof}
  Consider the double EGF of Theorem~\ref{tautheorem}, $a(y,z) = \frac{1}{1-z}exp(\sum_{i|k} 
  \frac{y^i-z^i}{i})$.
  Recall, the coefficient of $y^sz^t$ is the probability that $\pi^k \in S_{s+t}$ has $s$ fixed
  points.  Now, for any given $s$, we can find the probability that $\pi^k \in S_n$ has $s$
  fixed points (in the limit as $n \to \infty$), by evaluating $\lim_{z \to 1^-} (1-z)a(y,z)$.
  The result is the EGF $exp( \sum_{i|k} \frac{y^i - 1}{i} )$.
\end{proof}

However, the above requires us to have 256 terms inside of the exponentiation,
for there are 256 positive integers dividing 1,081,080, and we will need to know
the coefficient of the $c^{th}$ term for at least 1000 terms. Therefore, we are compelled
to leave this as a challenge for the computer algebra community.

Meanwhile, we performed the following experiment. We generated 10,000 random
permutations $\pi$ from $S_{10,000}$ and raised $\pi$ to the $k$th power for the
values of $k$ listed. Then we calculated $c$, the number of fixed points of $\pi^k$, and
determined if a search of the first 1/64th of the domain would reveal no fixed points.
That probability is given by 
$$ \left ( 1 - c/n \right )^{n/64} \approx e^{-c/64}$$
and taking the arithmetic mean over all experiments, one obtains

\begin{tabular}{rlll}
& No fixed points & One or more & \\
$k=1$ &  0.985041 & 0.014959 & Random \\
$k=1000000$ & 0.797284 & 0.202716 & Alice \\
$k=1081079$ & 0.984409 & 0.015591 & Bob\\
$k=1081080$ & 0.418335 & 0.581665 & Charlie\\
\end{tabular}

Perhaps this is unsurprising, as in the case of Charlie, we expect 256 fixed points,
and so it would be surprising if all of those were missing from a part of the domain
equal to 1/64th of the total domain in size. On the other hand, for Bob we expect
only 2 fixed points, and it is exceptional that we find one by accident. 

Finally, we observe that if there is an equal probability of an adversary being
presented with a random cipher from $S_{2^{64}}$ or 3-DES in the key of one of
our three users, iterated to their exponent, then the success probability of the
attacker would be for Alice 59.39\%, for Bob 50.03\%, and for Charlie 78.34\%.
Note in each case, we check only $2^{64}/64=2^{58}$ plaintexts, and so this attack is 
$2^{112}/2^{58}=2^{54}$ times faster than brute-force.

\paragraph{A General Maxim:} If a permutation must be iterated for some reason,
then it should be iterated a prime number of times, to avoid fixed points.

\subsection{A Key Recovery Attack}

Consider the cipher given by
$$ F_{k_1,k_2}(p)=E_{k_1}(E_{k_2}^{(n)}(E_{k_1}( p ))) = c$$
where $k_1$ and $k_2$ are keys, and $E$ is encryption
with a block
cipher (let $D_k(c)=p$ denote decryption). 
If $E$ is DES and $n=1$, then this is the 
``triple DES'' construction. Here, we consider that 
$E$ is AES-256 as an example, and $n$ is Charlie's number,
$1081080$. Then $F$ is a block cipher with 512-bit key and
128-bit plaintext block. We will refer to $k_1$ as the
outer key, and $k_2$ as the inner key.

Suppose an attacker had an oracle for $F$ that 
correctly encrypts with the correct $k_1$ and $k_2$
that the target is using. Call this oracle $\phi(p)$.
Observe that
$G_{k_3}(x) = D_{k_3}( \phi( D_{k_3}( x ) ) )$
will have $G_{k_3}(x)=E_{k_2}^{(n)}(x)$ if and only
if $k_3=k_1$. Thus if we can correctly guess the outer key,
we have an oracle for the $n$th iteration of encryption
under the inner key. If $k_3\neq k_1$, then provided that
$E_{k_1}$ is computationally indistinguishable from a
random permutation from $S_{2^{128}}$ when $k_1$ is chosen uniformly at
random (a standard assumption) then $G_{k_3}(x)$ also
behaves as a random permutation.

Thus, for $k_1=k_3$, we can expect $G_{k_3}(x)$ to 
behave like Charlie's cipher in the previous section,
and for $k_1\neq k_3$, we can expect $G_{k_3}(x)$
to behave like a random permutation in the previous section.

Let one run of the distinguishing attack signify guessing
all possible $k_3$ values, and executing the previous
section's attack for each key. If ``random'' is indicated (i.e. no
fixed point found), then
we reject the $k_3$ but if ``real'' is indicated (i.e. at
least one fixed point found), then we add $k_3$ to a 
``candidate list.''

After one run of this distinguishing attack, we would 
have a candidate list of outer keys of expected size
$$(0.014959)(2^{256} - 1) + (0.581665)(1)$$
where the success probabilities are given in the
previous section, for the attack on Charlie.

If we repeat the distinguisher attack on these
candidate keys, taking care to use a distinct set of
plaintexts in our search, the success probabilities
will be the same. This non-overlapping property of the
plaintext search could be enforced by selecting the
six highest-order bits of the plaintext to be the value
of $n$. After $n$ runs, we would expect the list to 
contain
$$(0.014959)^n(2^{256} - 1) + (0.581665)^n(1)$$
candidate keys.

Of course, the true $k_3=k_1$ key will be present
with probability $0.581665^n$. Next, for each key $k_c$ on
the candidate list, we will check all possible $2^{256}$ 
values of $k_2$ (denoted $k_x$), via checking if
$$ p = \phi( D_{k_c}( D_{k_x}^{(n)}( D_{k_c}( p ) ) ) ) $$
which will be true if $k_x=k_2$ and $k_c=k_1$. This
check should be made for roughly 4--6 plaintexts, to 
ensure that the match is not a coincidence. This necessity
arises from the fact that the cipher has a 512-bit key
and 128-bit plaintext. We will be very conservative,
and select 6.

The number of encryptions required for the $n$ runs is
\begin{eqnarray*}
(1081080+4)(\frac{2^{128}}{64})
(2^{256} + (0.014959)(2^{256}) + 
(0.014959)^2(2^{256}) + \\
(0.014959)^3(2^{256}) + 
\cdots + (0.014959)^n(2^{256})) \\
=(1081080+2)(2^{378})\frac{1 - (0.014959)^{n+1}}{1 - 0.014959}\\
= 2^{398.06579\cdots} (1 - 0.014959^{n+1})
\end{eqnarray*}
and for the second stage
\begin{eqnarray*}
  (6)(2)(2+1081080)(2^{256})(0.014959^n)(2^{256})=(2^{535.6290\cdots})(0.014959^n)\\
   = 2^{535.6290-6.062842n}
\end{eqnarray*}
for a success probability of $(0.581665)^n$.

\sloppy
Using Maple, we find that $n=23$ is optimal, leaving a candidate list of
$2^{116.555\cdots}$ possible keys, and requiring $2^{398.41207\cdots}$ encryptions,
but with success prob\-abil\-ity 
$(0.581665)^{23}\approx 2^{-17.98001\cdots}$. 
A brute-force search of the $2^{512}$ possible
keys would have $(6)(2)(1081082)2^{512}$ encryptions to perform, or $2^{535.629007\cdots}$.
Naturally, if a success probability of $2^{-17.98001\cdots}$ were desired, then only
$2^{517.649\cdots}$ encryptions would be needed for that brute-force search.

Therefore this attack is $2^{119.237}$ times faster than brute-force search.

\section{Conclusions}
In this paper, we presented a known theorem on the probabilities of random permutations having
given cycle structures and cycle counts, along with several useful corollaries.
To demonstrate the applicability of this technique to cryptanalysis, we have taken two
attacks which were heretofore presented at least partially heuristically, and made them
fully rigorous. It is hoped that other attacks which rely upon detecting these probabilities
via experimentation will be made rigorous as well, by calculation via EGFs and OGFs.
We also hope that we have demonstrated the utility of analytic combinatorics in general,
as well as EGFs and OGFs in particular.

We also presented a new attack, on very highly iterated permutations. While the scenario
is not reasonable, and it is only a distinguisher attack, it is also interesting that the
$\tau$ function occurs here. If a permutation should be highly iterated, it should be
iterated a prime number of times. However, the choice of 25 on the part of UNIX
designers was not bad, as $\tau(25)=3$. We also extended this to a key-recovery attack,
in an unusual context.  It is unclear in what situations such large
numbers of iterations would occur, but from a pure mathematical point of view,
the additional security granted by prime iteration counts is interesting.



\appendix

\section{Of Pure Mathematical Interest}

The authors encountered the following interesting connections with some
concepts in number theory, but they turned out to be not needed in the body
of the paper. We present them here for purely scholarly interest.

\subsection{The Sigma Divisor Function}

\begin{lemma}\label{sigmalemma}
The sum
$\sum_{i | k} 1/i = \frac{1}{k} \sigma(k)$ where 
both $i$ and $k$ are positive integers, and where
$\sigma(k)$ is the divisor function (i.e. the sum
of the positive integers which divide $k$).
\end{lemma}

\begin{proof}
$$ \sum_{i | k} 1/i 
= \frac{k}{k} \sum_{i | k} 1/i 
= \frac{1}{k} \sum_{i | k} k/i
= \frac{1}{k} \sum_{i | k} i
= \frac{1}{k} \sigma(k)$$
\end{proof}

\begin{corollary}
Let $\pi$ be a permutation taken at random from $S_n$.
The probability that $\pi^k$ is a derangement is $e^{-\sigma(k)/k}$, in 
the limit as $n\rightarrow \infty$.
\end{corollary}

\begin{proof}
Let $D$ be the set of positive integers dividing $k$.
From Corollary~\ref{divisorcorollary}, we know that $x$ is a fixed point of
$\pi^k$ if and only if $x$ is in a cycle of length found
in $D$ for $\pi$.

We will use Corollary~\ref{prohibcycles}, with $\overline{A}=D$.
We obtain the probability is $ e^{-\sum_{i\in D} 1/i}$,

and Lemma~\ref{sigmalemma} gives the desired result.
\end{proof}

Note that substituting $A=\set{1}$ into the above
yields the same result as Corollary~\ref{derangedcorollary}.

\subsection{Ap\'ery's Constant}

Corollary~\ref{PropC5} provides an amusing connection with Riemann's zeta
function.  Recall, for complex $s$, the infinite series, $\sum_{n \geq 1} 1/n^s$
defines the ``zeta function'' $\zeta(s)$, provided the series converges.

\begin{corollary}
  The probability that a random permutation (in the limit as the size grows
  to infinity) does not contain cycles of square length is:
  \[
    e^{-\sum_{i\geq 1} 1/i^2} = e^{-\zeta(2)} = e^{-\pi^2/6} \approx   
    0.19302529,
  \]
or roughly $1/5$.
\end{corollary}

\begin{corollary}
  The probability that a random permutation (in the limit as the size grows
  to infinity) does not contain cycles of cube length is:
  $e^{-\zeta(3)} \approx 0.30057532$ 
\end{corollary}

Note, $\zeta(3)$ is known as Ap\'ery's Constant \cite{apery}, and occurs in certain quantum
electrodynamical calculations, but is better known to mathematicians as being
the probability that any three integers chosen at random will have no common factor
dividing them all \cite{wolframapery}.


\ignore{
\section{The Weak Cycle Count Structure Theorem}
In the body of the paper we have the Strong Cycle Structure Theorem, and make a special
case of $B=\Z^+$ and $A$ left arbitrary. It seems unfair not to do the same with the other
parameter. Therefore, in the special case that $A = \Z^+$, and $B$ is any subset of $\Z^+$, we obtain:

\begin{theorem}
The Weak Cycle Count Structure Theorem:

The combinatorial class $\fancy{P}^{(\Z^+, B)}$ has associated EGF, $P^{(\Z^+, B)}_e(z) = \beta\left(
\log\left(\frac{1}{1-z}\right)\right)$, that is,
$$\frac{1}{1-z} -  \left [ \sum_{k\not\in B} \frac{1}{k!}\left ( \log \left ( \frac{1}{1-z}\right ) \right )^k \right ]$$
\end{theorem}

\begin{proof} Using the EGF for $\Z^+$ and The Strong Cycle Structure Theorem, we obtain
\begin{eqnarray*}
  P^{(\Z^+, B)}_e(z) & = & \sum_{k \in B} \frac{1}{k!}\left(\log\left(\frac{1}{1-z}\right)\right)^k \\
  & = & \left [ \sum_{k=1}^{\infty} \frac{1}{k!}\left ( \log \left ( \frac{1}{1-z}\right ) \right )^k \right ]
  -  \left [ \sum_{k\not\in B} \frac{1}{k!}\left ( \log \left ( \frac{1}{1-z}\right ) \right )^k \right ]\\
  & = & \left [ exp(  \log \left ( \frac{1}{1-z}\right ) \right ]  -  \left [ \sum_{k\not\in B} \frac{1}{k!}\left ( \log \left ( \frac{1}{1-z}\right ) \right )^k \right ]\\
  & = & \frac{1}{1-z} -  \left [ \sum_{k\not\in B} \frac{1}{k!}\left ( \log \left ( \frac{1}{1-z}\right ) \right )^k \right ]
\end{eqnarray*}
\end{proof}

If we consider $B=\set{k}$ then the OGF is $z^k$ and the EGF is $z^k/k!$. Therefore it is not
unexpected that we obtain an EGF of 
$$\frac{1}{k!} \log \left ( \frac{1}{1-z} \right )^k$$ and thus expanding this as a proper power
series in $z$ allows us to calculate the probability that a permutation $\pi$ from $S_n$
consists of precisely $k$ orbits. It is the coefficient of the $k$th term of the expanded power
series EGF.

\section{Probable Garbage}
The cipher triple-DES has a 64-bit plaintext block and a 112-bit key \cite{tripleDES}.
Suppose some paranoid person decides to encrypt their messages by encrypting
under triple-DES $2^{17}$ times over. Then, 
$D_k(p)$ encrypts the plaintext $p$ under key $k$, and the person is using
$D_k^{131,072}(p)$.

Because $2^{17}$ is ``highly composite'', we will show that the $2^{17}$th
iteration of a permutation has many more small cycles than expected.
We will define a ``small cycle'' as one that is length $s$ or shorter.

\subsubsection{Formal Derivation}
Let $\pi$ be a permutation chosen at random from $S_n$, and all probabilities
here are stated in the limit as $n$ goes to infinity.
Pick a random $x$ in the domain of $\pi$. Of the cycles in the disjoint cycle notation
for of $\pi$, it is in the domain of exactly one of them, call it $\psi$. Denote by $\ell$,
the length of the cycle.
Then rewrite $\ell=2^ab$, 
where $a\ge 0$ and $2$ does not divide $b$. Of course,
if $\ell$ is odd, then $a=0$. In any case, we know 
$$\gcd(131,073, \ell)=\gcd(2^{17}, 2^ab)=2^{\min(17,a)}$$

This allows for the following possible outcomes
using Corollary~\ref{cyclesplittingcorollary}.
\begin{enumerate}
\item If $a=0$, then $\psi$ is replaced by another cycle of length $\ell=m$.
\item If $1 \le a \le 17$ then $\psi$ is replaced by $2^a$ cycles of length $b$.
\item If $a \ge 18$ then $\psi$ is replaced by $2^{17}$ cycles of length $2^{-17}\ell$.
\end{enumerate}

For the case $a=0$, then $x$ is in a cycle of length $s$ or shorter if and only
if $b=\ell \le s$. Since all cycle lengths under $\pi$
are equally likely (from Corollary~\ref{independentcorollary}), this occurs with probability
$s/n$, given $a=0$.

For $a=i$ with $1 \le i \le 17$, then $x$ is in a cycle of length $s$ or shorter if and only
if $b \le s$ or $\ell \le 2^{i}s$. Since all cycle lengths are equally likely, this occurs 
with probability $2^is/n$, given $a=i$.

For $a\ge 18$, then $x$ is in a cycle of length $s$ or shorter if and only if $b \le s$
or $\ell \le 2^{17}s$. Since all cycle lengths are equally likely, this occurs with
probability $2^{17}s/n$, given $a\ge 18$. Note, here and in the previous case
we have tacitly assumed $s < 2^{-17}n$ otherwise it is more proper to say
$\ell \le \min ( 2^{17}s, n )$ and so the probability of $\ell$ obeying that inequality
is more obviously $\min(2^{17}s, n)/n$.

So given any value of $a$ we know the probabilities, and so we can use conditional
probability formulas, because the values of $a$ are mutually exclusive. Let the event
$S$ be that $x$ is in a short cycle under $\pi^{131,073}$. Then
\begin{eqnarray*}
Pr\{S\} & = & \sum_{i=0}^{i=\floor{\log_2 n}} Pr\{S | a=i\}Pr\{a=i\} \\
& = & 
\left [ \left ( s/n \right ) \left ( 1/2 \right ) \right ] + 
\left [ \left ( \sum_{i=1}^{i=17} 2^{i}s/n \right ) \left ( 2^{-i-1} \right ) \right ] 
+ \left [ \left ( 2^{17}s/n \right ) \left ( 2^{-18} \right ) \right ]
= \frac{19s}{2n}
\end{eqnarray*}
where the probabilities for $a$ taking on the various values are known
from considering $\ell \mod 2^{18}$. 

For any given value of $a$, it is clear that the cycle length under $\pi^{131,072}$
will be uniformly distributed between $1$ and $s$. Therefore, it is true in general,
regardless of the probability distribution of $a$.

Suppose one operation consists of picking a random point, and checking its
cycle length, aborting early if $\ell > s$. Then if $x$ is in a short cycle, this will
take expected time $s/2$ evaluations of $\pi$. 
If $x$ is not in a short cycle, this will take $s$ evaluations of $\pi$. Thus the
duration of this operation is
$$ \frac{19s}{2n}\frac{s}{2} + \left ( 1-\frac{19s}{2n} \right ) s = s - 19s^2/4n$$

Now suppose we are willing to perform $w$ evaluations of $\pi$. Then
the number of iterations of this operation possible would be 
$\floor{ w(s -19s^2/4n)^{-1}}$, and the expected number of short cycles
found would be 
$$\floor{ w(s -19s^2/4n)^{-1}}\frac{19s}{2n}$$

shit, I'm just confused.
}
\ignore{
Now we will calculate the expected length of a short cycle. If $a=0$ then
we obtain 1 cycle of length uniformly distributed between $1$ and $s$.
Again, recall that for $\pi$, all the cycle lengths are equally likely by
Corollary~\ref{independentcorollary}.

For $a=i$ with $1\le i \le 17$, then we obtain $2^{i}$ cycles of length
uniformly distributed between $1$ and $s$. Likewise with $a\ge 18$,
we obtain $2^{17}$
}
\ignore{ 
\newpage

\section{The Old Abstract}
Random function and random permutation statistics
are very frequently used in cryptology.
However, it seems that for iterated permutations,
such as in slide attacks on block ciphers with
periodic structure (e.g. KeeLoq \cite{KeeLoq}),
the statistics needed in cryptography cannot
be found in the literature.

In this paper we first show
that the popular folklore method
for deriving a very basic result is incorrect
and we propose a simple fix for it,
that works however under a cryptographic assumption only.
Then we give a purely mathematical and
rigourous method to obtain this result
and several more general results,
by using the methods of analytic combinatorics.
By similar method we derive a number of
other known and new results, relevant to
an attack on the KeeLoq block cipher.

Very frequently, in cryptology, this kind of probabilities
are computed by doing simulations.
However since simulations on the real permutation
are totally infeasible due to its size,
an attacker in cryptology will run simulations
on smaller permutations.
But it appears that this type of approximation
was never studied and is likely to be mishandled.
Two questions need to be answered:
how precise is the result after running the simulations,
and how large should be the `toy''`random permutation to give a
good precision result in short time.
In this paper we show that, contrary to intuition,
random permutation statistics can be computed with
perfect accuracy already for \emph{very} small permutations
such as on $2^{5}$ elements or less,
and the optimal precision is achieved for smaller, not for bigger examples,
because of the sampling error that does \emph{not} depend of the permutation size.
Moreover, theory can give us exact results with precision
that infinitely superior to any computer experiment that can be handled.

This is very interesting,
because we will show on a contrived example
a cryptographic attack in which the mathematical theory
can give the answer,
but experimentation cannot if we do not have sufficient computing power.
\ignore{In other words a mathematic proof helps the attacker}

Old Keywords: block ciphers,
unbalanced Feistel ciphers,
random permutation statistics,
slide attacks,
KeeLoq.

\section{Old Introduction}

We denote the set of all positive integers by $\integer^+$.
In this section we study general permutations
$\pi:\{0,\ldots,N-1\}\to \{0,\ldots,N-1\}$
with $N$ elements,
for any $N\in \integer^+$, i.e. $N$ does not have
to be of the form $2^n$ (but in KeeLoq $N=2^{32}$).

Random permutations are frequently studied.
Nevertheless, we need to be careful, because
a number of facts about random permutations are counter-intuitive,
see for example Proposition \ref{PropC9} below.
Other intuitively obvious facts are not all obvious to prove rigourously,
as illustrated in the previous section when trying to prove
that the Proposition \ref{Prop1} holds for random permutations.
For these reasons, in this paper we will be
'marching on the shoulders of giants' and
will use some advanced tools,
in order to be able to compute
the exact answers to all our questions directly,
without fear of making a mistake,
and with a truly remarkable numerical precision.
\ignore{
}
We will be using methods of the modern analytic combinatorics.
It is a highly structured and very elegant theory,
that stems from the method of generating series attributed to Pierre-Simon Laplace,
and is greatly developed in the recent monograph book
by Flajolet and Sedgewick \cite{Flajolet}.
\ignore{
}

\section{Old Conclusion}

KeeLoq is a block cipher that is
in widespread use throughout the automobile industry
and is used by millions of people every day.
KeeLoq is a weak and simple cipher, and has several vulnerabilities.
Its greatest weakness is clearly the periodic property of KeeLoq
that allows many sliding attacks which are very efficient
\cite{BogdanovKeeLoq,BogdanovKeeLoq2,BogdanovKeeLoq3,KeeLoq,KeeLoqTatra,KeeLoqDunkelman}.
Their overall complexity simply does not depend on the number of rounds of the cipher.
Thus, efficient attacks on the full 528-round KeeLoq become possible.

In addition, KeeLoq has a very small block size of 32 bits which makes it
feasible for an attacker to
know and store more or less the whole code-book.
But when the attacker is given the whole code book of
$2^{32}$ plaintexts-ciphertext pairs,
KeeLoq becomes extremely weak.
Courtois Bard and Wagner showed that
it can be broken in time which is faster than
the time needed to compute the code-book
by successive application of the cipher.
This is the Slide-Determine attack from \cite{KeeLoq}
which is also the fastest key recovery attack currently known on KeeLoq.
In this paper we have described and analysed
an Improved Slide-Determine attack that is faster both
on average and in specific sub-cases (cf. Table  \ref{ComparisonTableKeeLoq2} below).
For example, for 30 $\%$ of all keys,
we can recover the key of the full cipher
with complexity equivalent to $2^{27}$ KeeLoq encryptions.
This attack does not work at all for a proportion of 37 $\%$
so called strong keys \cite{KeeLoq,KeeLoqTatra}.

In Table \ref{ComparisonTableKeeLoq2} we compare our result to
previously published attacks on KeeLoq
\cite{BogdanovKeeLoq,BogdanovKeeLoq2,BogdanovKeeLoq3,KeeLoq,KeeLoqTatra,KeeLoqDunkelman},
We can observe that the fastest attacks do unhappily require the
knowledge of the entire code-book.
Then, there are attacks that require $2^{16}$ known plaintexts
but are slower, and finally, there is brute force attacks that
are the slowest. In practice however, since they require only 2 known
plaintexts and are actually feasible,
brute force will be maybe the only attack
that will be executed in practice by hackers.
Knowing which is the fastest attack on one specific cipher,
and whether one can really break into cars and how,
should be secondary questions in a scientific paper.
Instead, in cryptanalysis we need to study a variety of attacks
on a variety of ciphers.

\ignore{
\vskip-6pt
\vskip-6pt
\begin{table}
  \caption{Comparison of our attacks to other attacks reported on KeeLoq.}
    \label{ComparisonTableKeeLoq2}
\vskip-7pt
\vskip-7pt
\vskip-7pt
$$
\begin{array}{|c|c|c|c|c|c|c|c|c|c|}
\hline
\mbox{Type of attack}    &\mbox{Data} &\mbox{Time} &\mbox{Memory} &\mbox{Reference} \\
\hline
\hline
\mbox{Brute Force}    &2~ \mbox{KP}  &2^{63}  &\mbox{small}  & \mbox{}\\
\hline
\hline
\mbox{Slide-Algebraic}    &2^{16} \mbox{KP} &\mathbf{2^{53}}  &\mbox{small}  & \mbox{Slide-Algebraic Attack 2 of \cite{KeeLoq}}\\
\hline
\mbox{Slide-Meet-in-the-Middle}    &2^{16} \mbox{KP} &\mathbf{2^{46}}  &\mbox{small}  & \mbox{Biham, Dunkelman {\em et al}} \cite{KeeLoqDunkelman}\\
\hline
\mbox{Slide-Meet-in-the-Middle}    &2^{16} \mbox{CP} &\mathbf{2^{45}}  &\mbox{small}  & \mbox{Biham, Dunkelman {\em et al}} \cite{KeeLoqDunkelman}\\
\hline
\hline
\mbox{Slide-Correlation}    &2^{32} \mbox{KP}  &2^{51}  &\mbox{16 Gb}  & \mbox{Bogdanov} \cite{BogdanovKeeLoq,BogdanovKeeLoq2}\\
\hline
\mbox{Slide-Fixed Points}    &2^{32} \mbox{KP} &2^{40}  
&\mbox{18 Gb}  & \mbox{Attack 4 in eprint/2007/062/}\\
\hline
\mbox{Slide-Cycle-Algebraic}    &2^{32} \mbox{KP} &2^{40}  &\mbox{18 Gb}  & \mbox{Attack A in \cite{KeeLoqTatra}, cf. also \cite{KeeLoq0}}\\
\hline
\mbox{Slide-Cycle-Correlation}    &2^{32} \mbox{KP} &2^{40}  &\mbox{18 Gb}  & \mbox{Attack B in \cite{KeeLoqTatra}, cf. also \cite{BogdanovKeeLoq2}}\\
\hline
\multicolumn{5}{c}{~} \\[-5pt]
\cline{5-5}
\multicolumn{4}{c|}{~} & \mbox{Two previous versions in \cite{KeeLoq}:}\\
\hline
\mbox{Slide-Determine}    &2^{32} \mbox{KP} &2^{31}  &\mbox{16 Gb}  & \mbox{A: for 63 $\%$ of all keys}\\
\hline
\mbox{Slide-Determine}    &2^{32} \mbox{KP} &2^{28}  &\mbox{16 Gb}  & \mbox{B: for 30 $\%$ of all keys}\\
\hline
\multicolumn{4}{c|}{~} & \mbox{Three new improved versions:}\\
\hline
\mbox{Slide-Determine}    &2^{32} \mbox{KP} &\mathbf{2^{30}}  &\mbox{16 Gb}  & \mbox{Average time, 63 $\%$ of keys}\\
\hline
\mbox{Slide-Determine}    &2^{32} \mbox{KP} &\mathbf{2^{27}}  &\mbox{16 Gb}  & \mbox{Realistic, 30 $\%$ of keys}\\
\hline
\mbox{Slide-Determine}    &2^{32} \mbox{KP} &\mathbf{2^{23}}  &\mbox{16 Gb}  & \mbox{Optimistic, 15 $\%$ of keys}\\
\hline
\end{array}
$$
\vskip-6pt
\vskip-6pt
\end{table}
\vskip-6pt
\vskip-6pt
Legend: The unit of time complexity here is one KeeLoq encryption.
\vskip6pt

\vfill

\vskip-7pt
\vskip-7pt
}

\begin{corollary}\label{CoroC21}
The expected number of fixed points
for the 512 rounds $f_k^8$ of KeeLoq
is exactly $\tau(8)=4$.
\end{corollary}

\vskip-6pt
\vskip-6pt
\begin{theorem}\label{PropC22}
Let $\pi$ be a random permutation and $k>0$.
The probability that $\pi^k$ has no fixed points
is:
\vskip-4pt
\vskip-4pt
$$
e^{-\sum_{i|k}
\frac{1}{i}}
~~~~\mbox{when}~~
N\to \infty
$$
\end{theorem}
\vskip-4pt

\noindent\emph{Justification:}
We combine Proposition \ref{PropC5} and Proposition \ref{PropC20}.

\vskip-6pt
\vskip-6pt
\begin{theorem}\label{PropC23}
Let $\pi$ be a random permutation and $j,k\in \integer++$.
The probability that $\pi^k$ has exactly $j$ fixed points
is:
\vskip-3pt
\vskip-4pt
\vskip-4pt
$$
e^{-\sum_{i|k}
\frac{1}{i}}
\cdot
S(j)
~~~~\mbox{when}~~
N\to \infty
~~~~\mbox{where}~~
S(j)=
\left[t^j\right]
exp\left(
\sum_{i | k}
\frac{t^i}{i}\right)
$$
\vskip-6pt
\vskip-6pt
\end{theorem}
\vskip-4pt

\noindent\emph{Justification:}
Let ${\cal A}=\{i\in\integer++ \mbox{~such that~} i|k\}$.
Any permutation can be split in a unique way
as follows:
we have one permutation
$\pi_{\cal A}$ that has only cycles
of type $i\in {\cal A}$,
and another permutation
$\pi_{{\cal A}^c}$
that has no cycles of
length $i\in {\cal A}$.
This, following the terminology of \cite{Flajolet} 
can be seen as a ``labelled product''
of two classes of permutations.
Following Section II.2. in \cite{Flajolet}),
the corresponding EGF are simply multiplied.
We need to compute them first.

Following Proposition \ref{PropC20}
the number of fixed points for $\pi^k$
is exactly the size of $\pi_{\cal A}$.
which must be equal to exactly $j$ for all the permutations
in the set we consider here.
We compute the EGF
of permutations $\pi_{\cal A}$
of size exactly $j$
and that have only cycles
of size $i\in {\cal A}$.
This EGF is simply obtained by taking
the part of degree exactly $j$
in the appropriate EGF written directly
following Proposition \ref{PropC3} as follows:

\vskip-5pt
\vskip-5pt
$$
g_{{\cal A},j}(t)=
z^j
\cdot
\left[t^j\right]
exp\left(
\sum_{i | k}
\frac{t^i}{i}\right)
$$
\vskip-3pt

Finally,
we multiply the result by the EGF from Proposition \ref{PropC5}:

\vskip-5pt
\vskip-5pt
$$
g(z)=
\frac{1}{1-z} \cdot
exp\left(-\sum_{i\in {\cal A}}
\frac{z^i}{i}\right)
\cdot
z^j
\cdot
\left[t^j\right]
exp\left(
\sum_{i | k}
\frac{t^i}{i}\right)
$$
\vskip-3pt

It remains to apply the Proposition \ref{PropC2} and we get the result.

}


                                                     %


\begin{thebibliography}{8}
\bibitem{apery}
Ap\'ery, R. ``Irrationalit\'e de $\zeta$(2) et $\zeta$(3)." Ast\'erisque 61, 11-13, 1979.

\bibitem{monograph} Gregory V. Bard.
  \emph{Algebraic Cryptanalysis}.
  Springer-Verlag. (Scheduled Release) 2009.

\bibitem{BardDis} Gregory V. Bard.
  \emph{Algorithms for the Solution of Linear and Polynomial
  Systems of Equations over Finite Fields, with Applications to Cryptanalysis}.
  PhD Dissertation. Department of Applied Mathematics and Scientific 
  Computation, University of Maryland at College Park. Defended April 30, 2007.

\bibitem{SATconv}
  Gregory V. Bard, Nicolas T. Courtois and Chris Jefferson:
  \sl Efficient Methods for Conversion and Solution of Sparse Systems of
  Low-Degree Multivariate Polynomials over GF(2) via SAT-Solvers, \rm
  Available at \url{http://eprint.iacr.org/2007/024/}.

\bibitem{KeeLoqDunkelman}
  Eli Biham, Orr Dunkelman, Sebastiaan Indesteege, Nathan Keller, Bart Preneel:
  How to Steal Cars --- A Practical Attack on KeeLoq,
  in Eurocrypt 2008, LNCS 4965, pp. 1-18, Springer, 2008.

\bibitem{BogdanovKeeLoq}
  Andrey Bogdanov: 
  \sl Cryptanalysis of the KeeLoq block cipher, \rm
  \url{http://eprint.iacr.org/2007/055}.

\bibitem{BogdanovKeeLoq2}
  Andrey Bogdanov:
  \sl Attacks on the KeeLoq Block Cipher and Authentication Systems, \rm
  3rd Conference on RFID Security 2007, RFIDSec 2007.

\bibitem{BogdanovKeeLoq3}
  Andrey Bogdanov:
  \sl Linear Slide Attacks on the KeeLoq Block Cipher, \rm
  The 3rd SKLOIS Conference on Information Security and Cryptology (Inscrypt 2007),
  LNCS, Springer-Verlag, 2007

\bibitem{KeeLoqTatra}
  Nicolas Courtois, Gregory V. Bard and Andrey Bogdanov:
  \sl Periodic Ciphers with Small Blocks
  and Cryptanalysis of KeeLoq, \rm
  In Tatra Mountains Mathematic Publications,
  post-proceedings of Tatracrypt 2007 conference, to appear in 2008.

\bibitem{KeeLoq0}
  Nicolas Courtois, Gregory V. Bard, David Wagner:
  \sl Algebraic and Slide Attacks on KeeLoq, \rm
  Older preprint with an incorrect specification of KeeLoq,
  available at \url{eprint.iacr.org/2007/062/}.

\bibitem{KeeLoq}
  Nicolas Courtois, Gregory V. Bard, David Wagner:
  \sl Algebraic and Slide Attacks on KeeLoq, \rm
  In FSE 2008, pp. 97-115, LNCS 5086, Springer.
  Older (partly out of date) preprint available at \url{eprint.iacr.org/2007/062/}.

\bibitem{Flajolet}
  Philippe Flajolet and Robert Sedgewick;
  \sl Analytic Combinatorics, \rm
  book of 807 pages, to apear in Cambridge University Press
  in the first half of 2008.
  Available {\bf in full} on the Internet,
  see \url{http://algo.inria.fr/flajolet/Publications/book.pdf}

\bibitem{coolbook} Simson Garfinkel and Gene Spafford. 
  \emph{Practical Unix \& Internet Security}.
  O'Reilly. 2nd edition. 1996.

\bibitem{Maple} Maple: A Computer Algebra System.
  \url{http://www.maplesoft.com/}.

\bibitem{MarkoRiedel} Marko R. Riedel,
  \sl Random Permutation Statistics, paper available on the internet, \rm
  at \url{http://www.geocities.com/markoriedelde/papers/randperms.pdf}.
  




































%
%
%


%






\bibitem{wolframapery}
Weisstein, Eric W. ``Ap\'ery's Constant.'' From MathWorld---A Wolfram Web Resource. 
\url{http://mathworld.wolfram.com/AperysConstant.html}
\end{thebibliography}
\end{document}